\documentclass[11pt,a4paper]{amsart}
\usepackage{amssymb}
\usepackage{hyperref}

\usepackage{graphicx,amssymb,amsfonts,epsfig,amsthm,a4,amsmath,url}
\usepackage[latin1]{inputenc}

\usepackage{amsmath, amsthm, amssymb, verbatim}
\usepackage{graphicx,amssymb,amsfonts,epsfig,amsthm,a4,amsmath,url,graphicx,enumerate}
\usepackage[latin1]{inputenc}
\usepackage{tikz}
\usepackage{textcomp}
\usetikzlibrary{matrix}
\usepackage{svg}

\newtheorem{thm}{Theorem}[section]
\newtheorem{cor}[thm]{Corollary}

\newtheorem{prop}[thm]{Proposition}
\theoremstyle{definition}

\theoremstyle{remark}

\numberwithin{equation}{section}

\begin{document}
\title[On Ahn-Hendrey-Kim-Oum question for twin-width of graphs with 6 vertices]{On Ahn-Hendrey-Kim-Oum question for twin-width of graphs with 6 vertices}
\author[Das]{Kajal Das}
\address{Indian Statistical Institute\\  203 Barrackpore Trunk Road\\Kolkata 700 108}
\email{kdas.math@gmail.com}

\maketitle
\textbf{Abstract:} Twin-width is a recently introduced graph parameter for finite graphs. It is an open problem to determine whether there is an $n$-vertex graph having twin-width at least $n/2$ (due to J. Ahn, K. Hendrey, D. Kim and S. Oum). In an earlier paper, the author showed that such a graph with less than equal to 5 vertices does not exist. In this article, we show that such a graph with 6 vertices does not exist. More precisely, we prove that 
each graph with 6 vertices has twin-width less than equal to 2.


\textbf{Key terms:}:  Finite graphs, Twin-width, Graphs with 6 vertices, Ahn-Hendrey-Kim-Oum conjecture.

\section{Introduction}

Twin-width is an invariant of graphs introduced in \cite{BKTW20}. It is defined for a finite simple graph, later it is extended for a simple infinite graph. It is used to study the parameterized complexity of graph algorithms. It has applications in logic, enumerative combinatorics etc. Recently, it has appeared in many articles (\cite{BGKTW21}, \cite{BGKTW21'}, \cite{BGMSTT21}, \cite{BKRT22}, \cite{BGTT22}, \cite{BCKKLT22}). Moreover, it has been studied in the context of finitely generated groups \cite{BGTT22}.   The computation of twin-width of a finite graph is extremely difficult. There are some results for some well known graphs,  for example, complete graphs, path graphs, cyclic graphs (or graphs with at most one cycle), Paley graphs, Caterpillar tree, planar graphs etc.

However, it is an open problem to determine whether there is an $n$-vertex graph having twin-width at least
$n/2$( (due to J. Ahn, K. Hendrey, D. Kim and S. Oum, see \cite{AHKO22}, page 3). In \cite{Das22}, the author proves the following theorem. 

\begin{thm}\cite{Das22}
Let $G$ be a graph with number of vertices less than equal to 5. Then, $G$  has twin-width less than equal to 2. In particular, the Ahn-Hendrey-Kim-Oum question is not true for graphs with vertices less than equal to 5. 
\end{thm}

In this article, we prove the following theorem. 

\begin{thm}(Main Theorem)\label{mainthm}
Let $G$ be a graph with number of vertices 6. Then, $G$  has twin-width less than equal to 2. In particular, the Ahn-Hendrey-Kim-Oum question is not true for graphs with number of vertices 6. 
\end{thm}

\subsection{Organization}
In Section 2, we introduce our necessary definitions, notations and abreviations. In Section 3, we discuss the ingredients for proving our main result. In Section 4, we prove our main theorem. 

\section{Preliminaries: some definitions, notations and abreviations:}

A \textit{trigraph} $G$ is a graph with a vertex set $V(G)$, a black edge set $E(G)$, and a red edge set $R(G)$ (the error
edges), where $E(G)$ and $R(G)$ are disjoint. The set of neighbours of a vertex $v$ in a
trigraph $G$, denoted by $N_G(v)$, consists of all the vertices adjacent to $v$ by a black or red edge. The degree of a vertex 
$v$ is defined by the number $\mid N_G(v)\mid$. A $d$-trigraph is a
trigraph $G$ such that the red graph $(V(G), R(G))$ has degree at most $d$. In this situation, we also
say that the trigraph has red degree at most $d$.

A \textit{contraction} or \textit{identification} in a
trigraph $G$ consists of merging two (non-necessarily adjacent) vertices $u$ and $v$ into a single
vertex $w$, and defining the edges of $G'$ (the new graph after contraction) in the following way: Every vertex of the symmetric
difference $N_G(u)\bigtriangleup N_G(v)$ is linked to $w$ by a red edge. Every vertex $x$ of the intersection $N_G(u)\cap N_G(v)$ is linked to $w$ by a black edge if both $ux\in E(G)$ and $vx\in E(G)$, and by a red edge otherwise. The rest of the edges (not incident to $u$ or $v$) remain unchanged.  Also, the vertices $u$ and $v$ (together with the edges incident to these vertices) are removed from the trigraph. 

A \textit{sequence of $d$-contractions} is a sequence of $d$-trigraphs $G_n, G_{n-1},\cdots, G_1$, where $G_n = G$, $G_1 = K_1$, where $K_1$ is the graph on a single vertex, and $G_{i-1}$ is obtained from $G_i$ by performing a single contraction of two (non-necessarily adjacent) vertices. We observe that $G_i$ has precisely $i$ vertices, for every $i\in \{1,\cdots, n\}$. The twin-width of $G$, denoted by $tww(G)$, is the minimum integer $d$ such that $G$ admits a $d$-sequence.

Now, we provide an example of a sequence contractions of a finite graph. In the sequence of graphs depicted below, we start with the given finite graph in the extreme left end and we label the vertices by $a, b, c, d, e, f, g$. The next diagram is the result of the contraction of the vertices $e$ and $f$ and in the resulting graph we label the new vertex by $ef$. In this way, we obtain a sequence of graphs by gradually contacting other vertices. The graph in the extreme left end of the second line of this sequence is obtained by contracting the vertices $ad$ and $g$ in the graph depicted in the extreme right end of the first line of the sequence.

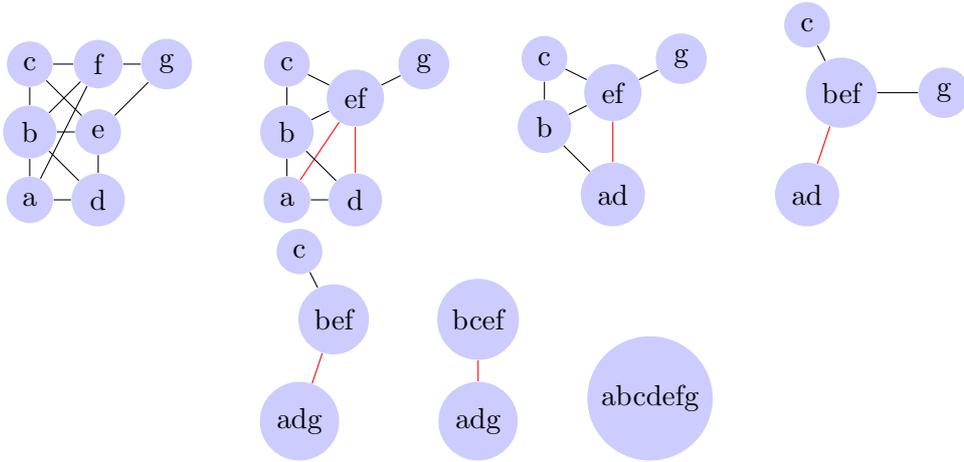
\begin{center}
\begin{tikzpicture}  
  [scale=.9,auto=center,every node/.style={circle,fill=blue!20}] 
    
  \node (a1) at (0,0) {a};  
  \node (a2) at (0,1)  {b};  
  \node (a3) at (0,2)  {c};  
  \node (a4) at (1,0) {d};  
  \node (a5) at (1,1)  {e};  
  \node (a6) at (1,2)  {f};  
  \node (a7) at (2,2)  {g};  
  
  \draw (a1) -- (a2); 
  \draw (a1) -- (a4);
  \draw (a1) -- (a6);
  \draw (a2) -- (a3); 
  \draw (a2) -- (a4); 
  \draw (a2) -- (a5);
  \draw (a2) -- (a6);
 \draw (a3) -- (a5);
 \draw (a3) -- (a6);
 \draw (a4) -- (a5);  
 \draw (a5) -- (a7);
 \draw (a6) -- (a7);  
\end{tikzpicture} 
\qquad
\begin{tikzpicture}  
  [scale=.9,auto=center,every node/.style={circle,fill=blue!20}] 
    
  \node (a1) at (0,0) {a};  
  \node (a2) at (0,1)  {b};  
  \node (a3) at (0,2)  {c};  
  \node (a4) at (1,0) {d};  
  \node (a5) at (1,1.5)  {ef};  
  \node (a7) at (2,2)  {g};  
  
  \draw (a1) -- (a2); 
  \draw (a1) -- (a4);
  \draw[red] (a1) -- (a5);
  \draw (a2) -- (a3); 
  \draw (a2) -- (a4); 
  \draw (a2) -- (a5);
  \draw (a3) -- (a5);

 \draw[red] (a4) -- (a5);  
 \draw (a5) -- (a7);

\end{tikzpicture}
\qquad
\begin{tikzpicture}  
  [scale=.9,auto=center,every node/.style={circle,fill=blue!20}] 
    
  \node (a1) at (1,0) {ad};  
  \node (a2) at (0,1)  {b};  
  \node (a3) at (0,2)  {c};  
 
  \node (a5) at (1,1.5)  {ef};  
  \node (a7) at (2,2)  {g};  
  
  \draw (a1) -- (a2); 
  \draw[red] (a1) -- (a5);
  \draw (a2) -- (a3); 

  \draw (a2) -- (a5);
  \draw (a3) -- (a5);

 \draw (a5) -- (a7);

\end{tikzpicture}
\qquad
\begin{tikzpicture}  
  [scale=.9,auto=center,every node/.style={circle,fill=blue!20}] 
    
  \node (a1) at (0,-0.5) {ad};  
  \node (a2) at (0.5,1)  {bef};  
  \node (a3) at (0,2)  {c};  
  \node (a7) at (2,1)  {g};  
  
  \draw[red] (a1) -- (a2); 
  \draw (a2) -- (a3); 
  \draw (a2) -- (a7);

\end{tikzpicture}
\qquad
\begin{tikzpicture}  
  [scale=.9,auto=center,every node/.style={circle,fill=blue!20}] 
    
  \node (a1) at (0,-0.5) {adg};  
  \node (a2) at (0.5,1)  {bef};  
  \node (a3) at (0,2)  {c};

  \draw[red] (a1) -- (a2); 
  \draw (a2) -- (a3);

\end{tikzpicture}
\qquad
\begin{tikzpicture}  
  [scale=.9,auto=center,every node/.style={circle,fill=blue!20}] 
    
  \node (a1) at (0,-0.5) {adg};  
  \node (a2) at (0,1)  {bcef};

  \draw[red] (a1) -- (a2); 
\end{tikzpicture}
\qquad
\begin{tikzpicture}  
  [scale=.4,auto=center,every node/.style={circle,fill=blue!20}] 
    
  \node (a1) at (0,0) {abcdefg};

\end{tikzpicture}

\end{center}
Moreover, we will draw every contraction sequence by this fashion in this article. We end this section 
by defining twin-width of an infinite graph. It is defined by the maximum of the twin-widths of its induced finite subgraphs. 

\section{Some ingredients for proving our result}

In this section, we discuss the twin-width of some well known graphs which will be useful for proving our main result. 

\subsection{Ingredient 1}

\begin{thm}\label{complete}
The complete graph with $n$-vertices, denoted by $K_n$,  and the complete bipartite with ${n,m}$-vertices, denoted by $K_{n,m}$,  have twin-width zero.
\end{thm}

\begin{prop}\label{ovcta}
If a graph has a vertex which is connected to every other vertices, then the graph has twin-width less than equal to 2.
\end{prop}



\subsection{Ingredient 2}

\begin{thm}\label{onecycle}\cite{AHKO22}
If every component of a graph $G$ has at most one cycle,
then $tww(G)\leq 2$.
\end{thm}

We obtain the following corollary from the above-mentioned theorem.
\begin{cor}\label{cyclicgraph}
The cyclic graph with $n$-vertices (denoted by $C_n$) has twin-width less than equal to 2.
\end{cor}

\subsection{Ingredient 3}

A \textit{caterpillar tree} is a tree in which all the vertices are within distance 1 of a central path.
We draw an example of a Caterpillar tree below. 

\begin{center}
\begin{tikzpicture} 
[scale=1.5,auto=center,every node/.style={circle,fill=blue!20}] 
\node (a1) at (0,0) {a}; 
\node (a2) at (1,0) {b}; 
\node (a3) at (2,0){c};
\node (a4) at (2,-1) {d};
\node (a5) at (3,0) {e};
\node (a6) at (2.5,1) {f};
\node (a7) at (3,1) {g};
\node (a8) at (3.5,1) {h};
\node (a9) at (2.5,-1) {i};
\node (a10) at (3.5,-1){j};
\node (a11) at (4,0) {k};

\draw (a1) -- (a2);
\draw (a2) -- (a3);
\draw (a3) -- (a4);
\draw (a3) -- (a5);
\draw (a5) -- (a6);
\draw (a5) -- (a7);
\draw (a5) -- (a8);
\draw (a5) -- (a9);
\draw (a5) -- (a10);
\draw (a5) -- (a11);
\end{tikzpicture}
\end{center}

\begin{thm}\label{caterpillar}\cite{AHKO22}
For a tree $T$, $tww(T)\leq 1$ if and only if $T$ is a Caterpillar tree. 
\end{thm}
The above-mentioned theorem gives rise to the following corollary. 
\begin{cor}\label{pathgraph}
The path graph with $n$ vertices, denoted by $P_n$, has twin-width less than equal to 1. 
\end{cor}








\subsection{Ingredient 4}

 \begin{thm}\label{4vertices}\cite{Das22}
The twin-width of a graph with 4 vertices is less than equal to 1.
\end{thm}

\begin{thm}\label{5vertices}\cite{Das22}
The twin-width of a graph with 5 vertices is less than equal to 2.
\end{thm}

\subsection{Ingredient 5}
 
 The \textit{complement of a graph} $G$ is a graph $H$ on the same vertices such that two distinct vertices of $H$ are adjacent if and only if they are not adjacent in $G$.  We obtain the following theorem from \cite{BKTW20} (see Subsection 4.1).
 
 \begin{thm}\label{complementgraph}
 Twin-width is invariant under complementation. 
 \end{thm}
 
\section{Twin-width of graphs with 6 vertices}

From \cite{CP84}, we have the list of graphs with 6 vertices. In this section, we show that each graph of this list has twin width less than equal to 2, i.e., we prove our main theorem \ref{mainthm}. 

1.
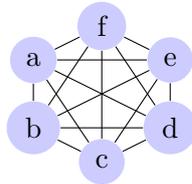
\begin{center}

\end{center} 
It is a complete graph. Therefore, the twin-width of this graph is zero. 

2.
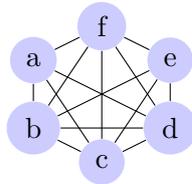
\begin{center}
%
\end{center} 
`b' (or `c' or `d' or `f') is connected with every other vertices. Therefore, the graph has twin-width less than equal to 2 by Proposition \ref{ovcta}. 

3.
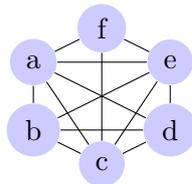
\begin{center}
%
\end{center} 
`a' (or `c' or `e') is connected with every other vertices. Therefore, the graph has twin-width less than equal to 2 by Proposition \ref{ovcta}. 

4.
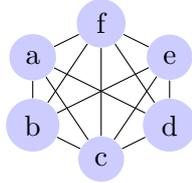
\begin{center}
%
\end{center} 
`c' (or `f') is connected with every other vertices. Therefore, the graph has twin-width less than equal to 2 by Proposition \ref{ovcta}. 

5.
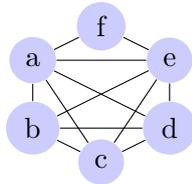
\begin{center}
%
\end{center}
`a' (or `e') is connected with every other vertices. Therefore, the graph has twin-width less than equal to 2 by Proposition \ref{ovcta}. 

6.
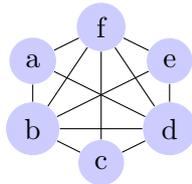
\begin{center}
%
\end{center} 
`b' (or `d' or `f') is connected with every other vertices. Therefore, the graph has twin-width less than equal to 2 by Proposition \ref{ovcta}. 

7.
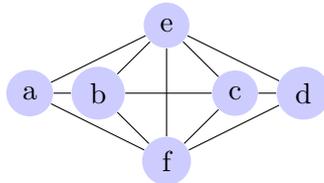
\begin{center}
%
\end{center} 
`e' (or `f') is connected to every other vertices. Therefore, the graph has twin-width less than equal to 2 by Proposition \ref{ovcta}. 

8.
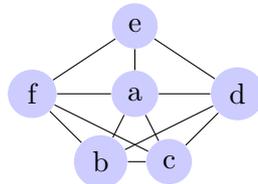
\begin{center}
%
\end{center} 
`a' is connected with every other vertices. Therefore, the graph has twin-width less than equal to 2 by Proposition \ref{ovcta}. 

9.
\begin{center}
%
\end{center} 
`d' is connected with every other vertices. Therefore, the graph has twin-width less than equal to 2 by Proposition \ref{ovcta} . 

11.
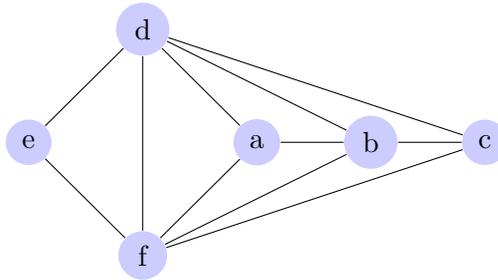
\begin{center}
%
\end{center} 
`d' (or `f') is connected with every other vertices. Therefore, the graph has twin-width less than equal to 2 by Proposition \ref{ovcta} . 

12.
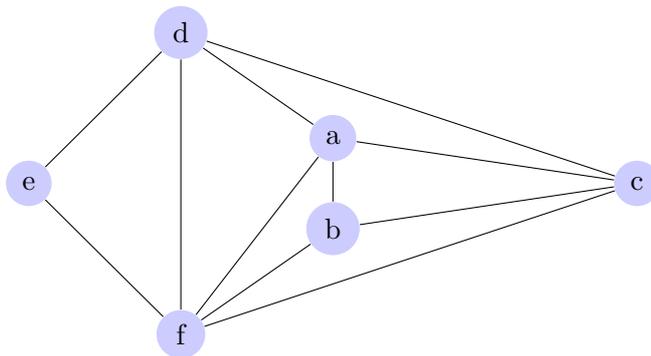
\begin{center}
%
\end{center} 
`f' is connected with every other vertices. Therefore, the graph has twin-width less than equal to 2 by Proposition \ref{ovcta} . 

13.
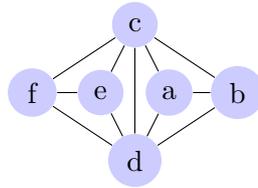
\begin{center}
%
\end{center} 
`c' (or `d') is connected with every other vertices. Therefore, the graph has twin-width less than equal to 2 by Proposition \ref{ovcta} . 

14.
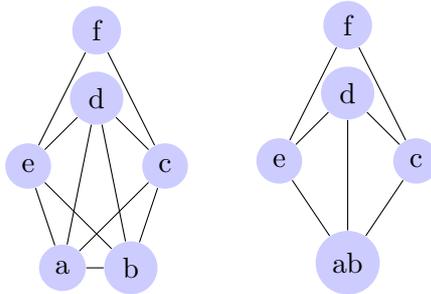
\begin{center}
%
\end{center}
After the contraction of the vertices `a' and `b', we obtain a graph with 5 vertices and with all black edges.  Therefore, the graph has twin-width less than equal to 2 by Theorem \ref{5vertices}.

15. 
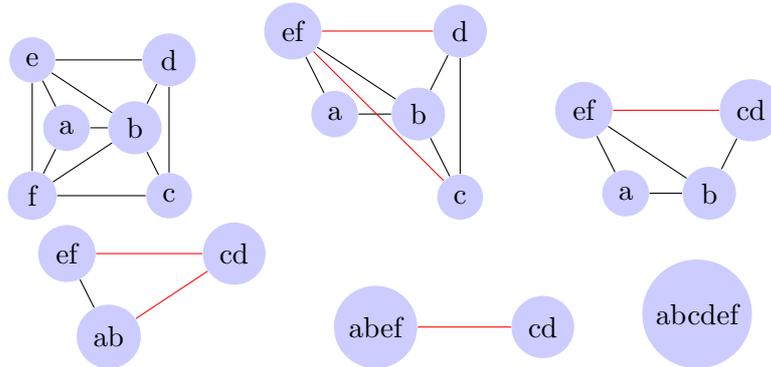
\begin{center}
%
\end{center} 
We obtain a 2-contraction sequence of the graph. Therefore, the graph has twin-width less than equal to 2.

16.
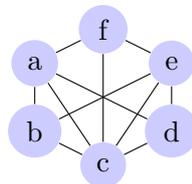
\begin{center}
%
\end{center} 
`c' is connected with every other vertices. Therefore, the graph has twin-width less than equal to 2 by Proposition \ref{ovcta}. 

17.
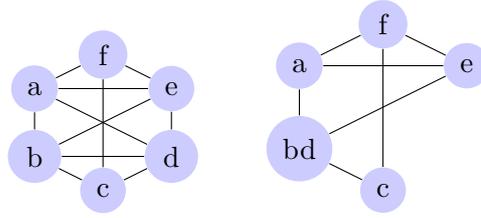
\begin{center}
%
\end{center} 
After the contraction of the vertices `b' and `d', we obtain a graph with 5 vertices and with all black edges. Therefore, the graph has twin-width less than equal to 2 by Theorem \ref{5vertices}.

18.
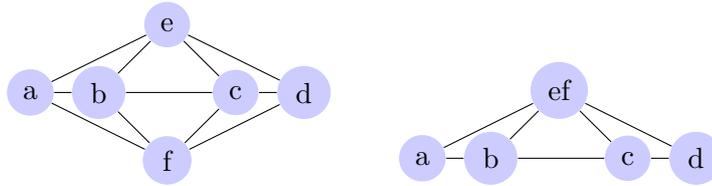
\begin{center}
%
\end{center} 
After the contraction of the vertices `e' and `f', we obtain a graph with 5 vertices and with all black edges. Therefore, the graph has twin-width less than equal to 2 by Theorem \ref{5vertices}.

19.
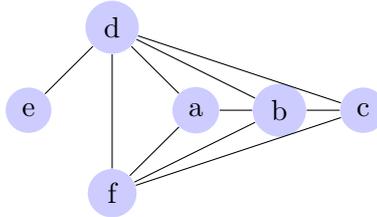
\begin{center}
%
\end{center} 
`d' is connected with every other vertices. Therefore, the graph has twin-width less than equal to 2 by Proposition \ref{ovcta}.

20.
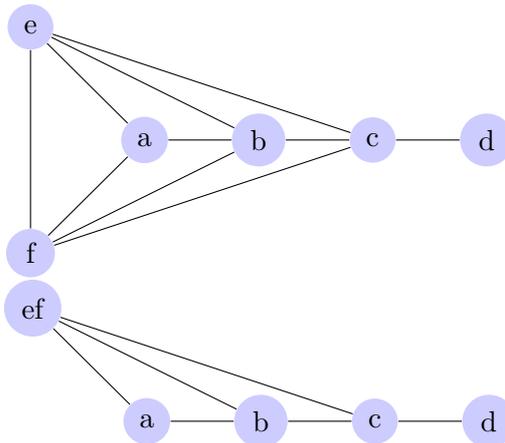
\begin{center}
%
\end{center} 
After the contraction of the vertices `e' and `f', we obtain a graph with 5 vertices and with all black edges. Therefore, the graph has twin-width less than equal to 2 by Theorem \ref{5vertices}.

21.
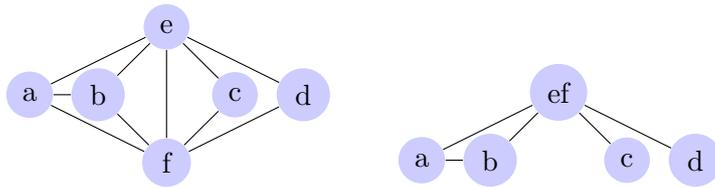
\begin{center}
%
\end{center} 
After the contraction of the vertices `e' and `f', we obtain a graph with 5 vertices and with all black edges. Therefore, the graph has twin-width less than equal to 2 by Theorem \ref{5vertices}.

22.
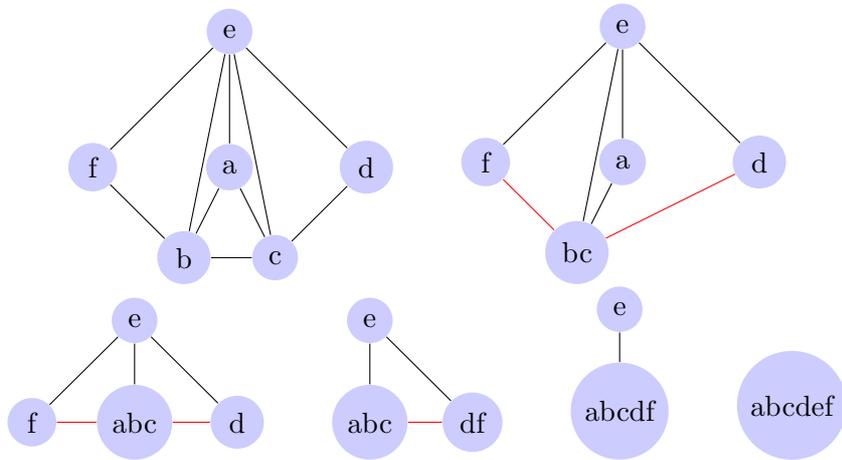
\begin{center}
%
\end{center} 
We obtain a 2-contraction sequence of the graph. Therefore, the graph has twin-width less than equal to 2. 

23.
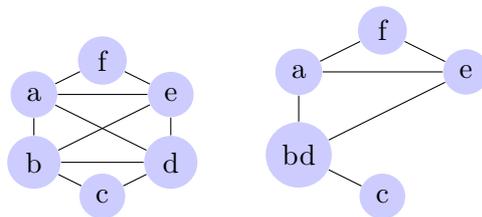
\begin{center}
%
\end{center}
After the contraction of the vertices `b' and `d', we obtain a graph with 5 vertices and with all black edges. Therefore, the graph has twin-width less than equal to 2 by Theorem \ref{5vertices}.
 
24.
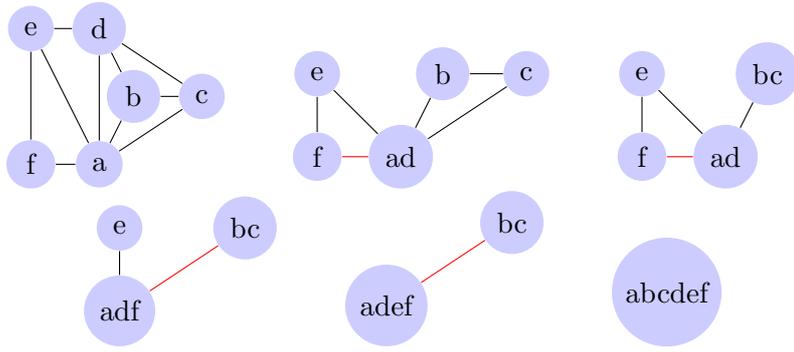
\begin{center}
%
\end{center} 
We obtain a 1-contraction sequence of the graph. Therefore, the graph has twin-width less than equal to 1.

25.
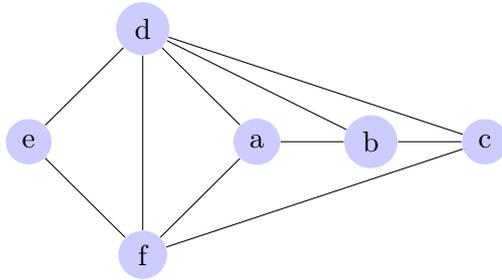
\begin{center}
%
\end{center} 
`d' is connected with every other vertices. Therefore, the graph has twin-width less than equal to 2 by Proposition \ref{ovcta}.

26.
\begin{center}
%
\end{center} 
`a' is connected with every other vertices. Therefore, the graph has twin-width less than equal to 2 by Proposition \ref{ovcta}. 

29.
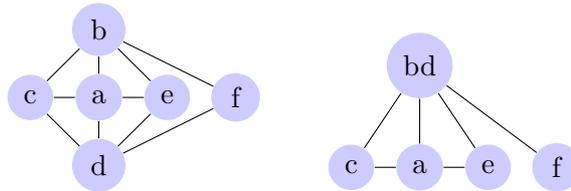
\begin{center}
%
\end{center}
After the contraction of the vertices `b' and `d', we obtain a graph with 5 vertices and with all black edges. Therefore, the graph has twin-width less than equal to 2 by Theorem \ref{5vertices}.

30.
\begin{center}
%
\end{center}
After the contraction of the vertices `d' and `f', we obtain a graph with 5 vertices and with all black edges. Therefore, the graph has twin-width less than equal to 2 by Theorem \ref{5vertices}.

35.
\begin{center}
%
\end{center}
After the contraction of the vertices `c' and `e', we obtain a graph with 5 vertices and with all black edges. Therefore, the graph has twin-width less than equal to 2 by Theorem \ref{5vertices}.

36.
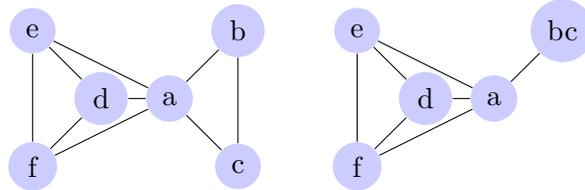
\begin{center}
%
\end{center}
After the contraction of the vertices `b' and `c', we obtain a graph with 5 vertices and with all black edges. Therefore, the graph has twin-width less than equal to 2 by Theorem \ref{5vertices}.

37.
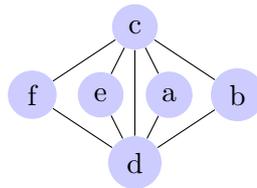
\begin{center}
%
\end{center} 
`c' is connected with every other vertices. Therefore, the graph has twin-width less than equal to 2 by Proposition \ref{ovcta}.

38.
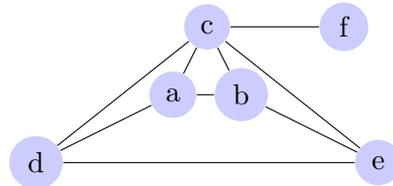
\begin{center}
%
\end{center} 
`c' is connected with every other vertices. Therefore, the graph has twin-width less than equal to 2 by Proposition \ref{ovcta}.

39.
\begin{center}
%
\end{center}
In the above picture, if we contract the vertices `d' and `f', we obtain a graph with 5 vertices and with all black edges. Therefore, the graph has twin-width less than equal to 2 by Theorem \ref{5vertices}.

46.
\begin{center}
%
\end{center} 
In the above picture, if we contract the vertices `a' and `e', we obtain a graph with 5 vertices and with all black edges. Therefore, the graph has twin-width less than equal to 2 by Theorem \ref{5vertices}.

50.
\begin{center}
%
\end{center}
In the above picture, if we contract the vertices `d' and `f', we obtain a graph with 5 vertices and with all black edges. Therefore, the graph has twin-width less than equal to 2 by Theorem \ref{5vertices}.

51.
\begin{center}
%
\end{center} 
In the above picture, if we contract the vertices `b' and `d', we obtain a graph with 5 vertices and with all black edges. Therefore, the graph has twin-width less than equal to 2 by Theorem \ref{5vertices}.

53.
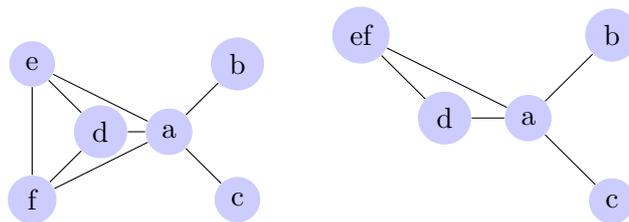
\begin{center}
%
\end{center}
In the above picture, if we contract the vertices `e' and `f', we obtain a graph with 5 vertices and with all black edges. Therefore, the graph has twin-width less than equal to 2 by Theorem \ref{5vertices}.

54. 
\begin{center}
%
\end{center}
In the above picture, if we contract the vertices `e' and `f', we obtain a graph with 5 vertices and with all black edges. Therefore, the graph has twin-width less than equal to 2 by Theorem \ref{5vertices}.

56.
\begin{center}
%
\end{center}
In the above picture, if we contract the vertices `b' and `d', we obtain a graph with 5 vertices and with all black edges. Therefore, the graph has twin-width less than equal to 2 by Theorem \ref{5vertices}.

58.
\begin{center}
%
\end{center}
In the above picture, if we contract the vertices `b' and `c', we obtain a graph with 5 vertices and with all black edges. Therefore, the graph has twin-width less than equal to 2 by Theorem \ref{5vertices}.

62.
\begin{center}
%
\end{center}
In the above picture, if we contract the vertices `b' and `c', we obtain a graph with 5 vertices and with all black edges. Therefore, the graph has twin-width less than equal to 2 by Theorem \ref{5vertices}.

63.
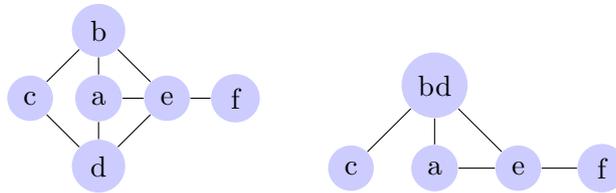
\begin{center}
%
\end{center}
In the above picture, if we contract the vertices `b' and `c', we obtain a graph with 5 vertices and with all black edges. Therefore, the graph has twin-width less than equal to 2 by Theorem \ref{5vertices}.

64.
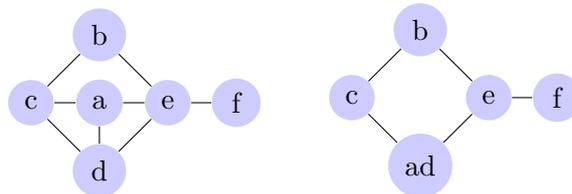
\begin{center}
%
\end{center}
In the above picture, if we contract the vertices `a' and `d', we obtain a graph with 5 vertices and with all black edges. Therefore, the graph has twin-width less than equal to 2 by Theorem \ref{5vertices}.

65.
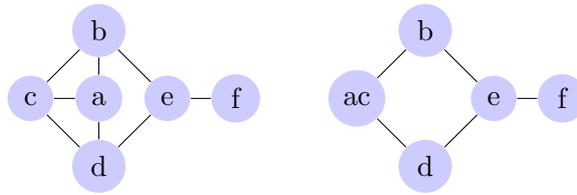
\begin{center}
%
\end{center}
In the above picture, if we contract the vertices `a' and `c', we obtain a graph with 5 vertices and with all black edges. Therefore, the graph has twin-width less than equal to 2 by Theorem \ref{5vertices}.

66.
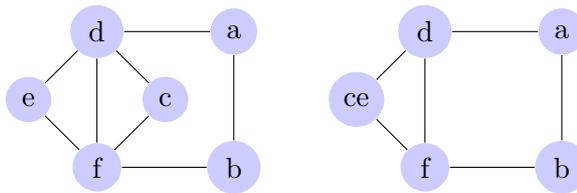
\begin{center}
%
\end{center}
In the above picture, if we contract the vertices `c' and `e', we obtain a graph with 5 vertices and with all black edges. Therefore, the graph has twin-width less than equal to 2 by Theorem \ref{5vertices}.

67.
\begin{center}
%
\end{center}
In the above picture, if we contract the vertices `c' and `e', we obtain a graph with 5 vertices and with all black edges. Therefore, the graph has twin-width less than equal to 2 by Theorem \ref{5vertices}.

70.
\begin{center}
%
\end{center}  
In the above picture, if we contract the vertices `a' and `e', we obtain a graph with 5 vertices and with all black edges. Therefore, the graph has twin-width less than equal to 2 by Theorem \ref{5vertices}.

72.
\begin{center}
%
\end{center} 
In the above picture, if we contract the vertices `c' and `d', we obtain a graph with 5 vertices and with all black edges. Therefore, the graph has twin-width less than equal to 2 by Theorem \ref{5vertices}.

74.
\begin{center}
%
\end{center}
In the above picture, if we contract the vertices `b' and `d', we obtain a graph with 5 vertices and with all black edges. Therefore, the graph has twin-width less than equal to 2 by Theorem \ref{5vertices}.

75.
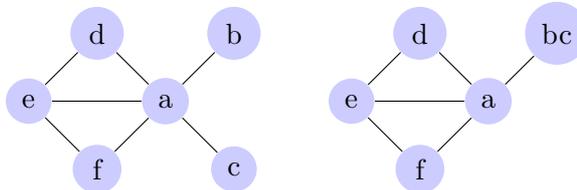
\begin{center}
%
\end{center}
In the above picture, if we contract the vertices `c' and `e', we obtain a graph with 5 vertices and with all black edges. Therefore, the graph has twin-width less than equal to 2 by Theorem \ref{5vertices}.

77.
\begin{center}
%
\end{center}
We obtain a 1-contraction sequence of the graph. The graph has twin-width less than equal to 1.

78.
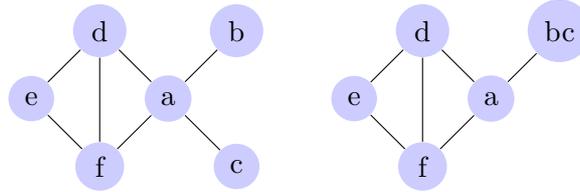
\begin{center}
%
\end{center}
In the above picture, if we contract the vertices `c' and `d', we obtain a graph with 5 vertices and with all black edges. Therefore, the graph has twin-width less than equal to 2 by Theorem \ref{5vertices}.

80.
\begin{center}
%
\end{center}
In the above picture, if we contract the vertices `a' and `c', we obtain a graph with 5 vertices and with all black edges. Therefore, the graph has twin-width less than equal to 2 by Theorem \ref{5vertices}.

81.
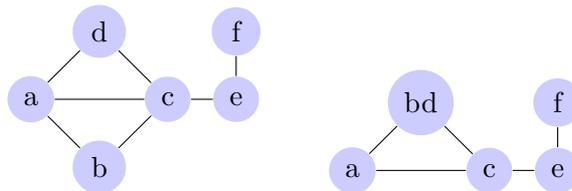
\begin{center}
%
\end{center}
In the above picture, if we contract the vertices `e' and `f', we obtain a graph with 5 vertices and with all black edges. Therefore, the graph has twin-width less than equal to 2 by Theorem \ref{5vertices}.

84.
\begin{center}
%
\end{center}
We obtain a 2-contraction sequence of the graph. Therefore, the graph has twin width less than equal to 2.

93.
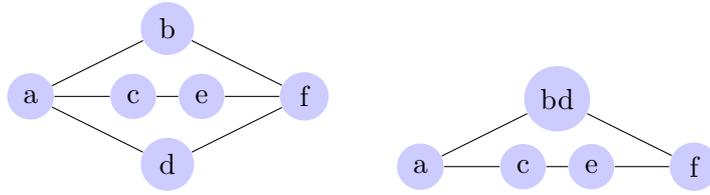
\begin{center}
%
\end{center}
After the contraction of the vertices `b' and `d', we obtain a graph with 5 vertices and with all black edges. Therefore, the graph has twin-width less than equal to 2 by Theorem \ref{5vertices}.

Now, we make a list of graphs which have at most one cycle. By Theorem \ref{onecycle}, these graphs have twin-width less than equal to 2.   

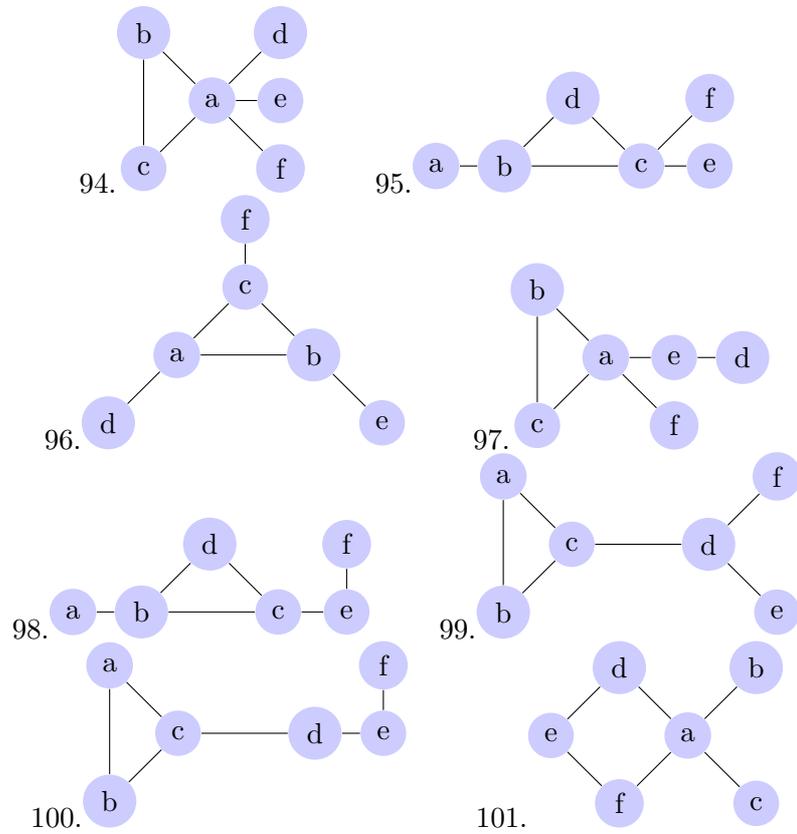
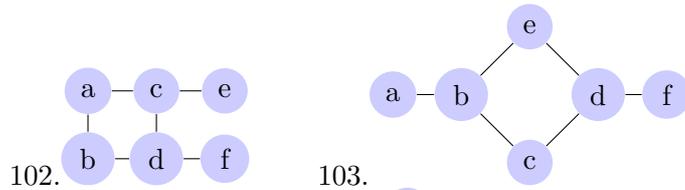
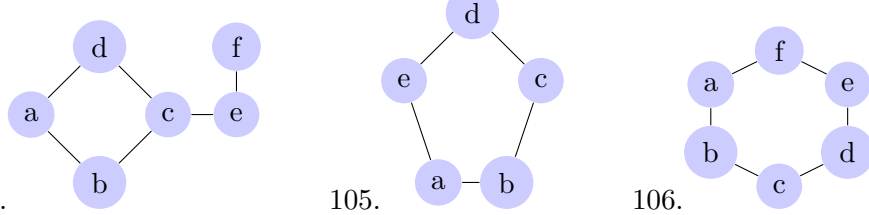
\begin{center}
94.%
\end{center}

We make a list of graphs which are Caterpillar tree. By Theorem \ref{caterpillar}, these graphs have twin-width less than equal to 1. 

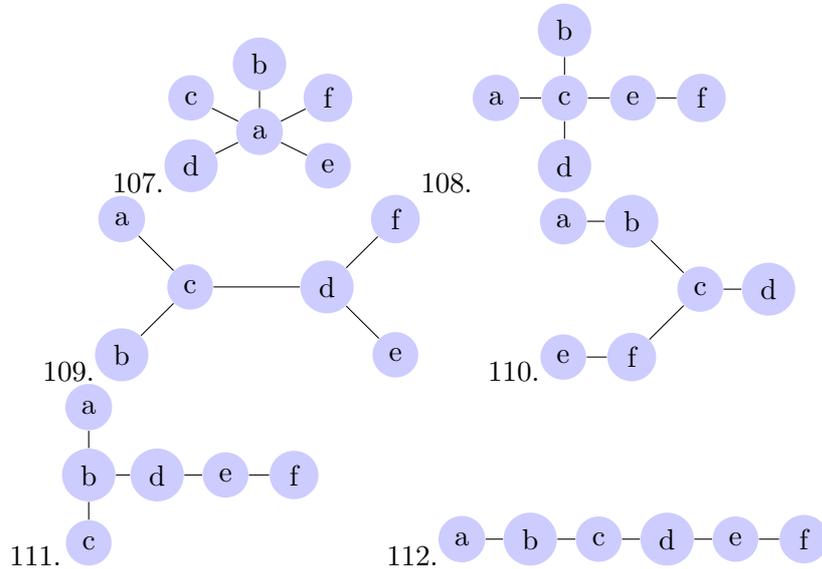
\begin{center}
107.%

\end{center}

\section{Acknowledgements}

I would like to thank National Board for Higher Mathematics (NBHM), India for supporting this project financially.

\end{document}